\documentclass[hidelinks,a4paper]{article}

\usepackage[T1]{fontenc}
\usepackage[utf8]{inputenc}
\usepackage{lmodern}
\usepackage{microtype}

\makeatletter

\usepackage[total={5.8in, 7.8in}, asymmetric, bindingoffset=0.4in]{geometry}
\usepackage{amsmath}
\usepackage{amssymb}
\usepackage{amsthm}
\usepackage{mathtools}
\usepackage{enumerate}
\usepackage[dvipsnames]{xcolor}

\usepackage{verbatim}

\usepackage{multirow}
\usepackage{empheq}
\usepackage{listings}
\usepackage{color}

\usepackage[square,comma,numbers]{natbib}
\usepackage{hyperref}
\usepackage{graphicx} 
\usepackage[font=small,labelfont=bf]{caption}
\numberwithin{equation}{section}
\usepackage{pgfplots}
\pgfplotsset{width=10cm,compat=1.9}

\usepackage[symbol]{footmisc}

\usepackage{fancyhdr}
\date{\vspace{-1em}\normalsize{\today}}

\def\cA{{\mathcal A}}
\def\cB{{\mathcal B}}
\def\cC{{\mathcal C}}

\def\cF{{\mathcal F}}

\def\cI{{\mathcal I}}
\def\cJ{{\mathcal J}}

\def\cL{{\mathcal L}}
\def\cM{{\mathcal M}}

\def\cO{{\mathcal O}}
\def\cP{{\mathcal P}}

\def\cT{{\mathcal T}}

\def\D{\mathbb{D}}
\def\E{\mathbb{E}}
\def\F{\mathbb{F}}

\def\H{\mathbb{H}}

\def\P{\mathbb{P}}
\def\R{\mathbb{R}}

\def\T{\mathbb{T}}
\def\Z{\mathbb{Z}}
\def\L{\mathbb{L}}

\def\eps{\epsilon}

\def\Om{\Omega}

\def\pmu{\partial_\mu}

\def\d{\mathrm{d}}

\def\balpha{{\boldsymbol{\alpha}}}

\def\hr{\widehat{\rho}}

\def\hc{\hat{c}}
\def\tc{\tilde{c}}

\def\bu {\bar{u} }

\def\ba{b^\alpha}

\def\sa{\sigma^\alpha}

\def\la{\ell^\alpha}
\def\lau{\ell^{\alpha_u}}

\theoremstyle{plain}
\newtheorem{theorem}{Theorem}[section]
\newtheorem{lemma}[theorem]{Lemma}
\newtheorem{corollary}[theorem]{Corollary}
\newtheorem{proposition}[theorem]{Proposition}
\newtheorem{assumption}[theorem]{Assumption}
\newtheorem{definition}[theorem]{Definition}

\newtheorem{remark}[theorem]{Remark}

\title{Viscosity Solutions for McKean-Vlasov Control on a torus\footnote{Partially supported 
by the National Science Foundation grant
 DMS 2106462.}}

\author{H. Mete Soner\footnote{Department of Operations Research and Financial
Engineering, Princeton University, Princeton, NJ, 08540, USA, email: 
{\tt soner@princeton.edu}. }
\and Qinxin Yan\footnote{Program in Applied and Computational
Mathematics, Princeton University, Princeton, NJ, 08540, USA, email: 
{\tt qy3953@princeton.edu}. }}

\date{\today}

\begin{document}
\maketitle
\abstract{
An optimal control problem in
the space of probability measures,  
and the viscosity solutions of the 
corresponding dynamic programming equations defined
using the intrinsic linear derivative are studied.
The value function is shown to be Lipschitz continuous
with respect to a novel smooth Fourier-Wasserstein metric.
A comparison result between the Lipschitz 
viscosity sub and super solutions
of the dynamic programming equation 
is proved using this metric,
characterizing the value function
as the unique Lipschitz viscosity solution. }
\vspace{3pt}

\vspace{3pt}


\section{Introduction}
\label{sec:intro}

McKean--Vlasov optimal control
is a part of the overarching program of
Lasry \& Lions \cite{LL1,LL2,LL} as articulated by Lions
through his College de France lectures \cite{Lions},
and independently initiated by Huang,
Malham\'{e}, \& Caines~\cite{HMRC}.
We refer the reader to
the classical
book of Carmona \& Delarue \cite{CD} and
to the lecture notes of Cardaliaguet \cite{C}
for detailed information and more references.

Main feature of the  McKean-Vlasov
type optimization  is the dependence 
of its evolution and cost not only on the position
of the state but also on its probability distribution, making the
set of probability measures as its state space.
Thus, the dynamic programming approach results
in nonlinear partial differential
equations set in the space of probability measures. 
Without common noise, 
they are first order Hamilton-Jacobi-Bellman
equations, and its 
Hamiltonian is defined only when 
derivative of the value function
is twice differentiable. In fact, this type
of unboundedness is almost
always the case for optimal control 
problems set in infinite dimensional
spaces \cite{GS} and is the main
new technical difficulty.  

These dynamic programming equations are analogous to the
coupled Hamilton-Jacobi and 
Fokker-Planck-Kolmogorov systems
that characterize the solutions of the 
mean-field games for which deep regularity
results are proved in  \cite{CDLL}
under some structural conditions.
However, in general the dynamic programming
equations for the McKean-Vlasov optimal control problems
are not expected to admit classical
solutions as shown in  subsection~\ref{ss:example} below, 
and a weak formulation is needed.

As the maximum principle is still the salient
feature in these settings as well, the
viscosity solutions of
Crandall \& Lions \cite{CL,CEL, CIL, FS}
is clearly the appropriate choice.  
However, due to the 
unboundedness of the Hamiltonian, 
original definition must
be modified. In fact, such modifications of 
viscosity solutions in infinite dimensional spaces
have already been studied extensively, and
the book \cite{GS}
provides an exhaustive account of these results.
Still, it is believed that
more can be achieved in the context
of McKean--Vlasov due to the special structure
of the set of probability measures.
Indeed, an approach developed by Lions
lifts the problems from the Wasserstein space
to a regular $\L^2$ space, and then exploits
the Hilbert structure to obtain new comparison results.  This procedure also 
delivers the novel \emph{Lions derivative}
which has many useful properties,
and we refer to  \cite{CD}
for its definition and more information.
This method is further developed
in several papers including \cite{BCFP,BCP,CGKPR1,PW1,PW}.
The choice of the appropriate notion
of a derivative is also explored in the recent paper \cite{GMS},
which then utilizes the deep connections to geometry
to prove uniqueness results for Hamiltonians
that are bounded in the sense discussed above.

Our main goals are to develop a viscosity theory 
directly on the space of probability measures using the linear derivative, provide
a comparison result, and obtain a
characterization of the value function
as the unique viscosity solution in a certain class
of functions.  A natural approach towards this goal 
is to project the problem onto finite-dimensional
spaces to leverage the already developed theory
on these structures. 
A second-order problem studied in \cite{CKLS}
provides a clear example of this approach
as its projections
exactly solve the projected finite dimensional equations.
However, in general these projections are only approximate solutions,
and \cite{CGKPR} uses the  Ekeland variational
principle together with 
Gaussian smoothed Wasserstein metrics as
gauge functions to control
the approximation errors.   
A different technical tool is developed in  \cite{BEZ}, and
\cite{GMS} studies the pure projection problem. 
Other approaches include the path-dependent 
 equations used in \cite{WZ}, gradient flows in \cite{CKT},
 convergence analysis in \cite{BeC} 
 and an optimal stopping problem in \cite{TTZ1,TTZ2}.
Recent paper \cite{CD1} exploits the semi-convexity, and also provides an extensive survey.

We on the other
hand employ the classical viscosity technique of doubling the variables
as done in \cite{MVJump} in lieu of projection.
The central difficulty of this approach is to
appropriately replace
the distance-square term $|x-y|^2$ used in the finite
dimensional comparison proofs with the square of a metric on
the space of measures. Thus, the crucial ingredient 
of our method is a novel Fourier-based
smooth metric  whose intriguing properties 
are studied in Section~\ref{sec:sdist}.
Our other main results are a comparison
between Lipschitz continuous sub and super
viscosity solutions, Theorem~\ref{th:compare} 
and the Lipschitz continuity 
of the value function with respect
to  a weaker metric, Theorem~\ref{th:Lipschitz}. 
Although the Lipschitz property
of the value function is  rather elementary 
for the Wasserstein metrics, it requires detailed
analysis for the Fourier based ones.  Indeed, a technical estimate,
Proposition~\ref{pro:Lipsde},
on the dependence of the solutions of the McKean--Vlasov
stochastic differential equation on the initial distribution
is needed for this property.

As our approach contains several new steps,
we study the simplest problem that allows us to showcase
its details and power concisely.
In particular, to ease the notation
we omit the dependence of all functions on the
time variable which can be added directly.  
Additionally, dynamics with jumps can
be included as done in \cite{MVJump}.
The compact structure of the 
torus is clearly a simplifying feature
as well.  In our accompanying paper \cite{SY}
we remove most of these restrictions
and study the extension of our method 
in higher dimensions.

The paper is organized as follow.  General 
structure and notations are given in the next section,
in Section~\ref{sec:problem} we define the problem
and state the assumptions. The main results are stated in
Section~\ref{sec:main}. 
We construct a family of Fourier-Wasserstein metrics
in Section~\ref{sec:sdist}. The comparison
result is proved in Section~\ref{sec:compare},
and the Lipschitz property in Section~\ref{sec:Lip}.
Standard results of dynamic programming
and viscosity property are proved in Section~\ref{sec:dpp} 
and respectively in Section~\ref{sec:viscosity}.

\section{Notations}
\label{sec:notations}
In this section, we summarize the notations and 
known results  used in the sequel.
We denote the dimension of the ambient space by $d$,
and  the finite horizon by $T>0$.
$\Z^d$ is the set of all $d$-tuples of integers.
$\T^d = \R^d /(2\pi\Z)^d$ is the $d$-dimensional torus
with the metric given by
$|x-y|_{\T^d}:=\inf_{k\in \Z^d}|x-y-2k\pi|$.
We use a filtered probability space 
$(\Omega,\F =({\cF}_t)_{0\leq t\leq T},\P)$ 
that supports Brownian motions.  We 
assume that initial filtration $\cF_0$ is rich
enough so that for any probability measure
on $\T^d$, there exists a random variable
on $\Om$ whose distribution is equal to this measure.

For a metric space $(E,d)$, 
$\cM(E)$ is the set of all Radon measures on $E$, and 
$\cP(E)$ denotes the set of all 
probability measures on $E$.
Let $\L^0(E)$ be the set of all $E$-valued
random variables. 
For $X \in \L^0(E)$,
$\cL (X)\in \cP(E)$ is the \emph{law of $X$}.

We denote the set of all continuous
real-valued functions on $E$ by $\cC(E)$, 
and the bounded ones by $\cC_b(E) \subset \cC(E)$.
We write $\cC(E,d)$ when the dependence on the 
metric is relevant, and $\cC(E \mapsto Y)$ if the range  
$Y$ is not the real numbers.
For a positive integer $n$, $\cC^n(E)$ is the set
of $n$-times continuously differentiable, real-valued
functions with the usual  norm $\|\cdot\|_{\cC^n}$
given by the sum of supremum norms of each derivative of 
order at most $n$.

We endow $\cM(E)$ with the 
 weak* topology $\sigma(\cP(E), \cC_b(E))$
and write $\mu_n\rightharpoonup \mu$, when
$\lim_{n \to \infty} \mu_n(f) =\mu(f)$ for every $f\in \cC_b(E)$.
Using the standard (linear) derivative
on the convex set $\cP(E)$, we say that
$\phi \in \cC(\cP(E))$ is 
\emph{continuously differentiable}
if there exists 
$\pmu \phi \in \cC(\cP(E) \mapsto \cC(E))$)
satisfying,
$$
\phi(\nu) 
=\phi(\mu) + \int_0^1\, \int_E\,
\pmu \phi(\mu + \tau (\nu-\mu))(x)\, (\nu-\mu)(\d x)\, 
\d \tau,
\qquad \forall\, \mu,\nu \in \cP(E).
$$

We set $\cO:=(0,T)\times\cP(\T^d)$.
For $\psi \in \cC(\overline{\cO})$ and $(t,\mu)\in \cO$, 
$\partial_t\psi(t,\mu)$ denotes the time derivative
evaluated at $(t,\mu)$, and  $\partial_\mu\psi(t,\mu)\in \cC(\T)$
denotes the derivative in the $\mu$-variable
again evaluated  at $(t,\mu)$. 
$\L^2(\T^d)$ is the set of
measurable  functions on $\T^d$
that are square integrable with respect to the Lebesgue measure,
with following orthonormal Fourier basis,
\begin{equation}
\label{eq:basis} 
e_k(x):= (2\pi)^{-\frac{d}{2}}\ e^{ik\cdot x},\qquad
x \in \T^d,\ k\in \Z^d,
\end{equation}
where $i =\sqrt{-1}$ and $z^*$ be the complex conjugate of $z$.
In particular, 
for any $\gamma\in \L^2(\T^d)$,
$$
\gamma= \sum_{k \in \Z^d}\ F_k(\gamma)\, e_k,
\quad
\text{where}
\quad
F_k(\gamma):= \int_{\T^d}\ \gamma(x) e_k^*(x)\ \d x, \ \
k \in \Z^d.
$$ 

Following metrics on $\cP(\T^d)$ are
given by their dual representations,
\begin{align*}
\rho_\lambda(\mu,\nu)&:= \sup\{ (\mu-\nu)(\psi)\ :\ \psi \in \H_\lambda(\T^d), \
\|\psi\|_\lambda \le 1\},\qquad \lambda \ge 1,\\
\hr_n(\mu,\nu)&:=  \sup\{ (\mu-\nu)(\psi)\ :\ \psi \in \cC^n(\T^d), \
\|\psi\|_{\cC^n} \le 1\},\qquad n=1,2,\ldots,
\end{align*}
where in view of Kantorovich duality, $\hr_1$ is the Wasserstein-one distance, and
for $\lambda \ge 1$, 
$$
\H_\lambda(\T^d):= \{ f\in \L^2(\T^d) : \|f\|_\lambda <\infty\},
\qquad
\|f\|_\lambda:= (\sum_{k \in \Z^d}\ (1+|k|^2)^\lambda\, |F_k(f)|^2)^\frac12.
$$
A Fourier representation of $\rho_\lambda$ is derived 
in Corollary~\ref{cor:rhoest}. 

It is well-known that
$\H_\lambda$ is the  classical Sobolev space with 
fractional derivatives.  Indeed, for any integer $n\ge 1$,
$\cC^n(\T^d) \subset \H_n(\T^d)=W^{n,2}(\T^d)$,
and $ \hr_n \le c_n \rho_n$ for some constant $c_n$. 
Moreover, by the embedding results, 
$\H_\lambda(\T^d) \subset C^n(\T^d)$ if 
$\lambda >n+\frac{d}{2}$. In particular,
we set 
\begin{equation}
\label{eq:star}
 n_*(d)=n_*:=3+ \lfloor\frac{d}{2}\rfloor,\quad
 \cC_*:= \cC^{n_*}(\T^d), \quad
\rho_*:=\rho_{n_*}, \quad
\hr_*:=\hr_{n_*},
\end{equation}
where $ \lfloor a \rfloor$
is the integer part of a real number $a$.
Then, $\H_{n_*}(\T^d) \subset \cC^2(\T^d)$. 

\section{McKean-Vlasov control }
\label{sec:problem}

In this section, we define 
the \emph{McKean-Vlasov optimal control} problem and
for a general  introduction,
we refer the reader to Chapter 6 in \cite{CD}.
Formally,  starting from $t \in [0,T]$,
the goal is to choose feedback controls $(\alpha_u(\cdot))_{u \in [t,T]}$
so as to minimize
$$
\int_t^T \E[ \ell(X_u,\cL(X_u), \alpha_u(X_u))]\ \d u
+ \varphi(\cL(X_T)),
$$
where $\ell$ is the 
\emph{running cost}, $\varphi$ is the
\emph{terminal cost}, $b,\sigma$ are given functions, and with a Brownian motion $B$,
$\d X_u = b(X_u,\cL(X_u), \alpha_u(X_u))\d u
+ \sigma(X_u,\cL(X_u), \alpha_u(X_u))\d B_u$.

We continue by defining this problem properly.

\subsection{Controlled processes}
\label{ss:cs}

Suppose that $A$ is a
closed Euclidean space and let
the \emph{control set} $\cC_a$ be
a subset of $\cC(\T^d \to A)$
containing all constant functions, and the
\emph{admissible controls} $\cA$ be the 
set of (deterministic) measurable 
functions $\balpha :[0,T] \mapsto \cC_a$.
We denote the value of any $ \balpha \in \cA$
at time $u \in [0,T]$
by  $\alpha_u \in \cC_a$. 
Given functions are the \emph{drift} vector $b=(b_1,\ldots,b_d) \in \R^d$,
the $d \times d'$ \emph{volatility} matrix $\sigma=(\sigma_{ij})$
with
$i=1,\ldots,d$,  $j=1,\ldots,d'$,
and the costs $\ell,\varphi$.
We continue by stating our standing regularity assumptions
on these functions,
$$
b_i, \sigma_{ij}, \ell : \T^d \times \cP(\T^d) \times A \mapsto  \R,
\qquad
\varphi :\cP(\T^d) \mapsto \R.
$$
Recall $\cC_*, \rho_*,\hr_*$ of \eqref{eq:star},
and for $\alpha \in \cC_a$, $x \in \T$, and $\mu \in \cP(\T^d)$,
set
$$
\ba(x,\mu):= b(x,\mu,\alpha(x)),\quad
\sa(x,\mu):= \sigma(x,\mu,\alpha(x)),\quad
\la(x,\mu):= \ell(x,\mu,\alpha(x)).
$$

\begin{assumption}[Regularity]
\label{a:regularity}
{\rm{There exists  $c_a < \infty$ such that 
for all $\alpha \in \cC_a$ and $\mu \in \cP(\T^d)$,
$$
\| \ba(\cdot,\mu)\|_{\cC_*}+
\| \sa(\cdot,\mu)\|_{\cC_*}+
\| \la(\cdot,\mu)\|_{\cC_*} \le c_a,
$$
and for $h=b,\sigma,\ell,\varphi$,
$$
|h(x,\mu,a)-h(x,\nu,a)| \le c_a\ \hr_*(\mu,\nu),
\qquad
\forall\, x \in \T^d,\, \mu,\nu \in \cP(\T^d),\,  a \in A.
$$
}}
\end{assumption}
\vspace{5pt}

Under this regularity condition,
for any $\balpha \in \cA$, $t \in [0,T]$,
and $\cF_t$ measurable, $\T^d$ valued random variable
 $\xi$ with $\mu =\cL(\xi)$, there is a 
 unique $\F$-adapted solution 
 $X_s^{t,\mu,\balpha}$ of the 
following \emph{McKean-Vlasov
stochastic differential equation},
\begin{equation}
\label{eq:mvsde}
X_s^{t,\mu,\balpha}
=\xi+\int_t^s
b^{\alpha_u}(X^{t,\mu,\balpha}_u, \cL^{t,\mu,\balpha}_u)\, \d u
+\int_t^s\sigma^{\alpha_u}(X^{t,\mu,\balpha}_u, \cL^{t,\mu,\balpha}_u)\, \d B_u,
\qquad s \in [t, T],
\end{equation}
where $\cL^{t,\mu,\balpha}_u =\cL(X^{t,\mu,\balpha}_u)$,
and $B$ is a $d'$  dimensional Brownian motion.

Although the solution
$X_u^{t,\mu,\balpha}$ depends on the choice of 
the initial condition $\xi$ and the Brownian increments
$(B_u-B_t)_{u \in [t,T]}$, 
as the Brownian increments are independent of $\cF_t$ and 
we consider feedback controls, the flow
 $(\cL_u^{t,\mu,\balpha})_{u\in[t,T]}$ depends only on the 
law $\mu = \cL(\xi)$ of the initial condition
 and not on $\xi$ itself.

 Clearly, the existence and uniqueness of 
 solutions of \eqref{eq:mvsde} can be 
 obtained under weaker assumptions. However,
 the stronger condition with $n_*$ 
 derivatives is needed for the comparison
 and the Lipschitz continuity results.
 We also emphasize that the regularity
 Assumption~\ref{a:regularity} 
 puts implicit regularity restrictions of the 
 control set $\cC_a$ as discussed 
in Remark~\ref{rem:justify} below.
 
\subsection{Problem}
\label{ss:problem}

Starting from $(t,\mu) \in \overline{\cO}$,
the  \emph{pay-off}
of a control process $\balpha \in \cA$
is given by,
\begin{equation}
\label{eq:j}
 J(t,\mu,\balpha):=\int_t^T
\E[\ell^{\alpha_u}(X^{t,\mu,\balpha}_u, \cL^{t,\mu,\balpha}_u)]\, \d u
+\varphi(\cL_T^{t,\mu,\balpha}),
\qquad \balpha \in \cA,\ (t,\mu) \in \overline{\cO}.
\end{equation}
Since
$\E[\ell^{\alpha_u}(X_u^{t,\mu,\balpha},\cL_u^{t,\mu,\balpha})]
=\cL_u^{t,\mu,\balpha}(\ell(\cdot,\cL_u^{t,\mu,\balpha},\alpha_u(\cdot)))$,
$J(t,\mu,\balpha)$ is 
a function of $\mu=\cL(\xi)$ independent of  the choice of  the initial
random variable $\xi$.
Although, this property, called \emph{law-invariance},
 holds directly
in our setting,
in general structures it is quite subtle.
We refer to
Proposition 2.4 of \cite{DPT}, and
Theorem 3.5 in  \cite{CGKPR1}
for its general proof,
and to Section 6.5 and Definition 6.27 of \cite{CD}
for a discussion. 

Then,
the McKean-Vlasov optimal
control problem  is to minimize the 
pay-off functional $J$ over $\balpha \in \cA$,
and the \emph{value function}
is given by,
$$
v(t,\mu):= \inf_{\balpha \in \cA}\ J(t,\mu,\balpha),
\qquad (t,\mu) \in \overline{\cO}.
$$

\begin{remark}
\label{rem:justify}
{\rm{
Suppose that 
$\cC_a = \{\alpha \in \cC_*(\T^d \to A) :\ \|\alpha\|_{\cC_*} \le c_0\}$
for some constant $c_0\ge0$.  Consider the class of functions
of the form
$h(x,\mu(f),a)$
for some $f \in \cC_*$,
and  $h:\T^d \times \R \times A \to \R$
satisfying
$\|h(\cdot,y,\cdot)\|_{\cC_*}+\|h(x,\cdot,a)\|_{1,\infty} \le c_1$
for every $x \in\T^d$, $y \in \R$, and $a \in A$,
for some $c_1\ge0$.
Then, $h^\alpha(x,\mu)= h(x,\mu(f),\alpha(x))$,
and $\|h^\alpha(\cdot,\mu)\|_{\cC_*}$
is less than a constant $c_a$ depending on $c_0,c_1$ and $n_*$.
Also, for every $x \in \T^d$,
$$
|h(x,\mu(f),\alpha(x))-h(x,\nu(f),\alpha(x))|\le
c_1 |(\mu-\nu)(f)| \le c_1 \|f\|_{\cC_*}\hr_*(\mu,\nu)\le
c_1c_0  \hr_*(\mu,\nu).
$$

Hence, this class of functions 
satisfy the regularity assumption.
More generally, under appropriate assumptions
functions  $h(x,\mu(f_1),\ldots,
\mu(f_m),a)$ 
with
$f_1,\ldots,f_m \in \cC_*(\T)$, 
and $h :\T^d \times \R^m \times A \to \R$ 
also satisfy the regularity assumption
with the above control set $\cC_a$. 
We emphasize that even when the coefficients
depend on $\mu$ only through $\mu(f_1),\ldots,\mu(f_m)$ of
the measure $\mu$, the value function in general
 is still infinite dimensional.

Assumptions made above hold
in a large class of examples 
studied in the mean-field games.  In particular,
for the Kuramoto problem studied in \cite{CCSo},
for some constants $\kappa, \sigma>0$,}}
$$
\ell(\mu,a)= \frac12 a^2 + \kappa[1- (\mu(\cos))^2- (\nu(\sin))^2],
\ \ b(x,\mu,a)=a,
\ \ \sigma(a)=\sigma.
$$
\end{remark}

\subsection{Dynamic programming principle} 
\label{ss:dpp}
We next state the dynamic programming principle
which is central to the viscosity approach to
optimal control.  A  general 
proof in a different setting is given in \cite{DPT}.
However, the continuity of the value
function proved in Section~\ref{sec:Lip}, and the standard techniques 
outlined in \cite{FS} 
allows for a simpler proof that we provide
in Section~\ref{sec:dpp}. 

\begin{theorem}[Dynamic programming]
\label{th:dpp}
For every $\mu \in \cP(\T^d)$
and $0 \le t \le \tau \le T$,
\begin{equation}
\label{eq:dpp}
v(t,\mu)=\inf_{\balpha\in\cA} \int_t^\tau 
 \E[\ell^{\alpha_u}(X^{t,\mu,\balpha}_u, \cL^{t,\mu,\balpha}_u)]\ \d u
+v(\tau, \cL_\tau^{t,\mu,\balpha}).
\end{equation}
\end{theorem}

It is well known that
the dynamic programming
can be used directly to show that the value function is a
viscosity solution of the 
dynamic programming equation
\begin{equation}
\label{eq:dpe}
- \partial_t v(t,\mu)=
H(\mu,\partial_\mu v(t,\mu)),
\qquad  t \in [0,T),\ \mu \in \cP(\T^d),
\end{equation}
where for $\gamma \in \cC^2(\T^d)$,
$\mu \in \cP(\T^d)$,  $x \in \T^d$ and $\alpha \in \cC_a$,
$$
H(\mu,\gamma) := \inf_{\alpha \in \cC_a} \big\{\mu 
\left(\ell^\alpha(\cdot,\mu)+ \cM^{\alpha,\mu}[\gamma](\cdot)
\right)\, \big\},
$$
$$
\cM^{\alpha,\mu}[\gamma](x) :=
b(x,\mu,\alpha(x))\cdot \partial_x\gamma(x)
+\sum_{i,j=1}^d  \sum_{l=1}^{d'}\ \sigma_{il}(x,\mu,\alpha(x))
\sigma_{jl}(x,\mu,\alpha(x)) 
\partial_{x_ix_j} \gamma(x).
$$
The value function
also trivially satisfies 
the following terminal condition,
\begin{equation}
\label{eq:initial}
v(T,\mu)=\varphi(\mu),\qquad \forall \ \mu \in \cP(\T^d).
\end{equation}

As the value function
is not necessarily differentiable,
a weak formulation is needed, and 
we use the notion of viscosity solutions.
The definition that we use 
 is exactly the classical one 
in which the auxiliary test functions
are continuously differentiable
functions on $\overline{\cO}=[0,T] \times \cP(\T^d)$, with the
linear derivative in $\cP(\T^d)$
recalled in Section~\ref{sec:notations}.
We continue by specifying the  auxiliary functions used in the definition of viscosity solutions.  

\begin{definition}
\label{def:smooth}
{\rm{We say that $\psi \in \cC(\overline{\cO})$ 
is a}} test function, {\rm{if
$\psi$ is continuously differentiable
with $\partial_\mu \psi(t,\mu) \in \cC^2(\T^d)$
for every $(t,\mu)\in \overline{\cO} $, and the map
$(t,\mu) \in \overline{\cO} \mapsto H(\mu,\partial_\mu\psi(t,\mu))$
is continuous. We denote  the set
of all test functions by  $\cC_s(\overline{\cO})$.}}
\end{definition}

\begin{definition}
\label{def:vis}
{\rm{A continuous function $u\in \cC(\overline{\cO})$
is a}} viscosity subsolution 
of \eqref{eq:dpe}, {\rm{if every $\psi \in \cC_s(\overline{\cO})$,
$(t_0,\mu_0) \in [0,T) \times \cP(\T^d)$, satisfying
$(u-\psi)(t_0,\mu_0) = \max_{\overline{\cO}}\, (u-\psi)$, 
also satisfies}}
$$
- \partial_t \psi(t_0,\mu_0) \le 
H(\mu_0,\partial_\mu \psi(t_0,\mu_0)).
$$

{\rm{A continuous function $w\in \cC(\overline{\cO})$
is a}} viscosity supersolution 
of \eqref{eq:dpe}, {\rm{if every $\psi \in \cC_s(\overline{\cO})$,
$(t_0,\mu_0) \in [0,T) \times \cP(\T^d)$, satisfying
$(w-\psi)(t_0,\mu_0) = \min_{\overline{\cO}}\, (w-\psi)$, 
also satisfies}}
$$
- \partial_t \psi(t_0,\mu_0) \ge 
H(\mu_0,\partial_\mu \psi(t_0,\mu_0)).
$$

{\rm{
Finally, $v\in \cC(\overline{\cO})$
is a}} viscosity solution of \eqref{eq:dpe}, 
{\rm{if it is both a sub and a
super solution.}}
\end{definition}

\section{Main results}
\label{sec:main}

Our main result is the 
characterization of the value
function as the unique continuous
viscosity solution of the dynamic programming equation 
\eqref{eq:dpe} and the terminal condition \eqref{eq:initial}. 

Recall the  metrics $\rho_*, \hr_*$ of \eqref{eq:star}.

\begin{theorem}[Comparison]
\label{th:compare}
Suppose that the regularity Assumption~\ref{a:regularity} holds, 
$u\in \cC(\overline{\cO})$ is a viscosity
subsolutionof \eqref{eq:dpe} and \eqref{eq:initial},
and $w \in \cC(\overline{\cO})$ is a viscosity
supersolution of \eqref{eq:dpe} and \eqref{eq:initial}.  If further
$u$ or $w$ is Lipschitz continuous
in the $\mu$-variable with respect
to the metric $\rho_*$, then $u \le v$ on $\overline{\cO}$.
\end{theorem}

Above comparison result is proved in Section~\ref{sec:compare} below.

\begin{theorem}[Continuity]
\label{th:Lipschitz}
Under the regularity Assumption~\ref{a:regularity},
there exists
a constant $L_v>0$ depending
only on the horizon $T$ and the constant 
$c_a$ of Assumption~\ref{a:regularity},  
so that
\begin{equation}
\label{eq:vLip}
|v(t,\mu)- v(s,\nu)| \le L_v \,\left[ \hr_*(\mu,\nu) + |t-s|^\frac12 \right],
\qquad \forall\, \mu,\nu \in \cP(\T^d), \, t,s \in [0,T].
\end{equation}
\end{theorem}
\noindent
This continuity result  proved in 
Section~\ref{sec:Lip} below,
also implies Lipschitz continuity
with respect to  $\rho_*$,
since  $\hr_* \le c_* \rho_*$ for some constant $c_*$.
The following result
follows directly from the standard viscosity
theory \cite{FS}, and its proof is given
in Section~\ref{sec:viscosity} below.

\begin{theorem}[Viscosity property]
\label{th:exist}
Under the regularity Assumption~\ref{a:regularity},
the value function is a 
viscosity solution of \eqref{eq:dpe} in $\cO$,
satisfying the terminal
condition \eqref{eq:initial}.

In particular, any continuous viscosity subsolution is less than or equal 
to the value function $v$, and any continuous viscosity supersolution is greater
than or equal to $v$.
\end{theorem}

\begin{remark}
\label{rem:lam}
{\rm{
In the comparison result,
we could use any metric $\rho_\lambda$ with $\lambda>2+\frac{d}{2}$.  
However,  our proof for Lipschitz
continuity requires us to employ the smaller metric $\hr_m$
and only for integer values of $m$.  This combination
of the results
dictates the global choice $\lambda =n_*$.}}
\end{remark}

\subsection{An example}
\label{ss:example}

In this subsection, we provide a simple example to 
illustrate the notation and also the need
for viscosity solutions.  We take $T=1$, $d=1$,
$A=\R$, $b(x,\mu,a)=a$, $\sigma\equiv 1$,
$\varphi \equiv 0$, and
$$
\ell(\mu,a):=\frac12 a^2+ L(m(\mu)),
\quad \text{where} \quad
m(\mu):= \int_\T x \ \mu(\d x),
$$
and $L:[-\pi,\pi] \to \R$ is a given Lipschitz function.
It can be shown that 
the value function of the above 
problem is independent of 
the control set  $\cC_a$, and  is given by,
$$
v(t,\mu) = w(t,m(\mu)),
\qquad (t,\mu) \in \overline{\cO},
$$
$$
w(t,y):= \inf_{\hat{\balpha} \in \widehat{\cA}}\ 
\hat{J}(t,y,\hat{\balpha}):=\inf_{\hat{\balpha}  \in \widehat{\cA}}
\int_t^1\ [\frac12(\hat{\alpha}_u)^2 + L(Y^{t,y,\hat{\balpha}}_u)]\ \d u,
\qquad (t,y) \in [0,1]\times \T,
$$
where $\widehat{\cA}$ is the set of
all measurable maps $\hat{\balpha}
: [0,1]\mapsto \T$, and
$ Y^{t,y,\hat{\balpha}}_u= y +\int_t^u  \hat{\alpha}_s \d s$.
It is well known that 
$w$ is the unique viscosity solution of
the Eikonal equation,
\begin{equation}
\label{eq:eikonal}
- \partial_t w(t,y) =- \frac12 (\partial_y w(t,y))^2+ L(y),\qquad y \in \T,
\end{equation}
and $w(1,\cdot)\equiv 0$.
Since $w$ is
not always differentiable,
we conclude that $v$ is not either,
and therefore  a weak theory is needed.
On the other hand, when $w$ is differentiable, we have
$$
\partial_\mu v(t,\mu)(x)= \partial_y w(t,m(\mu))\, x
\quad \Rightarrow \quad
\partial_x(\partial_\mu v(t,\mu)(x)) =\partial_y w(t,m(\mu)).
$$
Hence, by Jensen's inequality,
\begin{align*}
 H(\mu,\partial_\mu v(t,\mu)) &= \inf_{\alpha \in \cC_a}
 \int_\T \left(\frac12\alpha(x)^2 +\alpha(x) \partial_y w(t,m(\mu))\right)\mu(\d x)
 + L(m(\mu))\\
& \ge    \inf_{\alpha \in \cC_a}\left\{
\frac12 \left(\int_\T \alpha(x) \mu(\d x)\right)^2 +
 \left(\int_\T \alpha(x) \mu(\d x) \right) \partial_y w(t,m(\mu))\right\}
 + L(m(\mu))\\
 & =   \inf_{a \in \R}
\left\{\frac12 a^2 +
a \partial_y w(t,m(\mu))\right\}
 + L(m(\mu))\\
 &= -\frac12 ( \partial_y w(t,m(\mu)))^2 + L(m(\mu)).
 \end{align*}
As constant functions $\alpha \equiv a $
 are always in $\cC_a$,
 we also have the opposite inequality.
 Therefore, 
 $$
 H(\mu,\partial_\mu v(t,\mu)) =
-\frac12 ( \partial_y w(t,m(\mu))^2 + L(m(\mu)),
$$
for every $\cC_a$.
Since $\partial_tv(t,\mu)=\partial_t w(t,m(\mu))$,
the Eikonal equation \eqref{eq:eikonal} implies
that when $w$ is differentiable,
$v$ is a classical solution
of the dynamic programming equation
\eqref{eq:dpe}.

\section{Fourier-Wasserstein metrics}
\label{sec:sdist}

In this section, we study the properties of
the norms and the metric
$\rho_\lambda$ defined
in Section~\ref{sec:notations}.
Similar metrics are also
defined in \cite{Sobolev} using a dual representation
with Sobolev functions.

Recall that $z^*$ is the complex conjugate
of $z$, and the
orthonormal basis 
$\{e_k\}_{k \in \Z^d}$,
Fourier coefficients $F_k(f)$
are defined
in Section~\ref{sec:notations}.
For $\mu \in \cM(\T^d)$, $k \in \Z^d$, 
 we also set
$F_k(\mu):=  \mu(e^*_k)$.
As $\T^d$ is compact, $F_k(\mu)$ is finite for every $k$, and 
$F_0(\mu)=1$ for all $\mu \in \cP(\T^d)$.  

For $\lambda \ge 1$, we define a norm on $\cM({\T^d})$, dual to 
$\|\cdot\|_\lambda$ by,
$$
|\eta|_\lambda:= \sup\{ \eta(\psi)\ :\ \psi \in \H_\lambda({\T^d}), \
\|\psi\|_\lambda \le 1\},
\qquad \eta \in \cM({\T^d}),
$$
so that $\rho_\lambda(\mu,\nu)=  |\mu-\nu|_\lambda$.

\begin{lemma}
\label{lem:Hl}
For $\lambda > \frac{d}{2}$, $\eta \in \cM({\T^d})$, 
$|\eta|_\lambda <\infty$  and
has the following dual representation,
\begin{equation}
\label{eq:etal}
  |\eta|_\lambda= (\, \sum_{k \in \Z^d}\, (1+|k|^2)^{-\lambda}
   |F_k(\eta)|^2)^{\frac12}.
\end{equation}
\end{lemma}
\begin{proof}
We first note  as $ 2\lambda>d$,
$c_\lambda:=  \sum_{k \in \Z^d} (1+|k|^2)^{-\lambda}<\infty$.
Let $d(\eta)$ be the expression in
the right hand side of \eqref{eq:etal}
and $TV(\eta)$ be the total variation of the measure $\eta$.  
Then, $|F_k(\eta)| \le TV(\eta)$
and therefore, 
$d(\eta) \le \sum_{k \in \Z^d} (1+|k|^2)^{-\lambda} \, TV(\eta) 
=c_\lambda TV(\eta)$.

For $\psi \in \cC({\T^d})$, the Fourier representation 
$\psi=\sum_{k \in \Z^d} F_k(\psi)e_k$ implies that,
\begin{align}
\label{eq:1}
\eta(\psi)
&=\sum_{k\in \Z^d}\, F_k(\psi)\, \eta(e_k)
=\sum_{k\in \Z^d}\,  F_k(\psi)\, F_k^*(\eta)\\
\nonumber
&=  \sum_{k\in \Z^d}\,  [(1+|k|^2)^\frac{\lambda}{2}\, F_k(\psi)]\, 
[ (1+|k|^2)^{-\frac{\lambda}{2}}\, F_k^*(\eta)]\\
\nonumber
& \leq  ( \sum_{k\in \Z^d}\, (1+|k|^2)^\lambda {|F_k(\psi)|}^2)^\frac{1}{2}\
(\sum_{k\in \Z^d}\,  (1+|k|^2)^{-\lambda}\, |F_k^*(\eta)|^2 )^\frac12
= \|\psi\|_\lambda \,d(\eta).
\end{align}
In view of the definition of $|\cdot|_\lambda$,
$|\eta|_\lambda \le d(\eta)$, for any $\eta \in \cM({\T^d})$.

To prove the opposite inequality, fix $\eta \in \cM({\T^d})$
and  define a function $\tilde{\psi}$ by,
$$
\tilde{\psi}(x) := \sum_{k\in \Z^d}\, 
(1+|k|^2)^{-\lambda}\, F_{k}(\eta)\, e_k(x),
\quad
\Rightarrow
\quad
F_k(\tilde{\psi}) =(1+|k|^2)^{-\lambda}\, F_{k}(\eta),
\ \ k\in \Z^d.
$$
Since $c_\lambda<\infty$, $\tilde{\psi}$ is well-defined.
Moreover, 
$$
\|\tilde{\psi}\|_\lambda^2 
= \sum_{k \in \Z^d}
(1+|k|^2)^\lambda |F_k(\tilde{\psi})|^2
= \sum_{k \in \Z^d}
(1+|k|^2)^{-\lambda}\, |F_{k}(\eta)|^2
=d^2(\eta) <\infty.
$$
Hence, $\tilde{\psi} \in \H_\lambda({\T^d})$,
and by \eqref{eq:1},
$$
\eta(\tilde{\psi}) = 
\sum_{k\in \Z^d}\,  F_k(\tilde{\psi})\, F_{k}^*(\eta)
=\sum_{k\in \Z^d}\,  (1+|k|^2)^{-\lambda\, }\, | F_{k}(\eta)|^2
=d^2(\eta)
= \|\tilde{\psi}\|_\lambda \ d(\eta).
$$
As $\eta(\tilde{\psi})  \le 
|\eta|_{\lambda} \|\tilde{\psi}\|_\lambda$ by the definition of $|\cdot|_{\lambda}$, 
we have 
$d(\eta) \|\tilde{\psi}\|_\lambda = \eta(\tilde{\psi})  \le 
|\eta|_{\lambda} \|\tilde{\psi}\|_\lambda $.
\end{proof}

An immediate corollary is the following.

\begin{corollary}
\label{cor:rhoest}
For any $\lambda > \frac{d}{2}$,
$\rho_\lambda$ is a metric on $\cP(\T^d)$
with a dual representation,
$$
\rho_\lambda(\mu,\nu)= \max\{\ (\mu-\nu)(\psi)\ :\
\|\psi\|_\lambda\le 1\ \}=
(\, \sum_{k \in \Z^d}\, (1+|k|^2)^{-\lambda}\, |F_k(\mu-\nu)|^2\, )^{\frac12}.
$$
\end{corollary}
\begin{proof}
The dual representation follows directly from 
the previous lemma.
Suppose that
$\rho_\lambda(\mu,\nu)=0$, then $F_k(\mu)=F_k(\nu)$ 
for every $k \in \Z^d$.
As $\mu,\nu$ have the same Fourier series, we conclude that $\mu=\nu$. 
The fact that $\rho_\lambda$ is a metric now follows
from the dual representation.
\end{proof}

The following provides
a connection between the two metrics we consider.
Also with $m=1$, it implies that the classical
Wasserstein one metric $\hr_1$ is dominated by $\rho_1$.

\begin{lemma}
\label{lem:metric}
For any  integer $m\ge1$, there exists $c_{m,d}>0$,
such that
$\hr_m(\mu,\nu) \leq  c_{m,d}\, \rho_m(\mu,\nu)$
for every  $\mu,\nu \in \cP(\T^d)$.
\end{lemma}
\begin{proof} Fix the $m \ge1$ and let $D^m \psi$ be the $m$-th order
derivatives of $\psi \in  \cC^m(\T^d)$. 
Then, since $|k|^{2m} |F_k(\psi)|^2=|F_k(D^m \psi)|^2$,
$$
\sum_{k \in \Z^d} |k|^{2m} |F_k(\psi)|^2=
\sum_{k \in \Z^d}  |F_k(D^m \psi)|^2 =
\|\\D^m \psi\|^2_{\L^2(\T^d)}
\le d^m \ (2\pi)^d \|\psi\|^2_{\cC^m(\T^d)}.
$$
As 
$(1+|k|^2)^m \le 2^m(1+|k|^{2m})$,
for any $k\in \Z^d$, $\psi \in \cC^m(\T^d)$,
$$
\|\psi\|_m^2
\le 2^m  \sum_{k \in \Z^d} |F_k(\psi)|^2
+ 2^m \sum_{k \in \Z^d} |k|^{2m} |F_k(\psi)|^2
\le c_{m,d}^2\|\psi\|^2_{\cC^m(\T^d)},
$$
where $c_{m,d}^2=  2^m [1+  d^m  (2\pi)^d] $.  Hence,
\begin{align*}
\hr_m(\mu,\nu)&= \sup \{ (\mu-\nu)(\psi)\ :\  \|\psi\|_{\cC^m(\T^d)} \le 1\}\\
&\le  \sup \{ (\mu-\nu)(\psi)\ :\  \psi \in \cC^m(\T^d), \
\|\psi\|_m \le c_{m,d}\}\\
&\le  \sup \{ (\mu-\nu)(\psi)\ :\  \psi \in \H_m(\T^d), \
\|\psi\|_m \le c_{m,d}\} = c_{m,d} \rho_m(\mu,\nu).
\end{align*}
\end{proof}

Our next result is on the differentiability of $\rho_\lambda$.
Recall the test functions $\cC_s(\overline{\cO})$ of
Definition~\ref{def:smooth},  $n_*(d)$ of \eqref{eq:star},
and the basis $e_k$ of Section~\ref{sec:notations}.

\begin{lemma}
\label{lem:diff}
Fix $\lambda > \frac{d}{2}$,
$\nu \in \cP(\T^d)$ and set $h(\mu):= \frac12 \rho_\lambda^2(\mu,\nu)$.
Then, 
$$
\partial_\mu h(\mu)(x)
= \sum_{k\in \Z^d}\
(1+ |k|^2)^{-\lambda}\, F_k(\mu-\nu)
\, e^*_k(x),
\qquad x \in \T^d,
$$
and $\|\partial_\mu h(\mu)\|_\lambda = \rho_\lambda(\mu,\nu)$.
Moreover, if $\lambda = n_*(d)$, then
$\partial_\mu h(\mu) \in \cC^2(\T^d)$.
\end{lemma}
\begin{proof}
Fix  $\nu \in \cP(\T^d)$.
For each $k \in \Z^d$, set $a_k(\mu):=
\frac12|F_k(\mu-\nu)|^2$.
Then, we directly calculate that
$\partial_\mu a_k(\mu) (\cdot)
= F_k(\mu-\nu)\, e^*_k(\cdot)$.
Then, for any $x \in \T^d$,
$$
\partial_\mu h(\mu)(x) 
= \sum_{ k \in \Z^d} \,
(1+ |k|^2)^{-\lambda}\ \partial_\mu a_k(\mu) (x)
= \sum_{ k \in \Z^d} \,
(1+ |k|^2)^{-\lambda}\, F_k(\mu-\nu)\, e^*_k(x).
$$
The above formula implies that 
$F_k(\partial_\mu h(\mu))= (1+ |k|^2)^{-\lambda}F^*_k(\mu-\nu)$ 
for every $k \in \Z^d$.
Hence,
$$
\|\partial_\mu h(\mu)\|_\lambda^2 =
\sum_{k\in \Z^d}
(1+ |k|^2)^\lambda\, |F_k(\partial_\mu h(\mu))|^2=
\sum_{k\in \Z^d}
(1+ |k|^2)^{-\lambda}\, |F_k(\mu-\nu)|^2
= \rho_\lambda^2(\mu,\nu).
$$

In view of the Sobolev
embedding of $\H_{n_*}({\T^d})$ into $\cC^2({\T^d})$,
$\partial_\mu h(\mu) \in \cC^2(\T^d)$.
\end{proof}

\section{Comparison}
\label{sec:compare}

In this section we prove Theorem~\ref{th:compare}
in several steps.  Recall the test functions  
$\cC_s(\overline{\cO})$
of Definition~\ref{def:smooth},
and $n_*, \rho_*$ of \eqref{eq:star}.
Then, $2(n_*-2) \ge d+1$, and consequently, 
\begin{equation}
\label{eq:cd}
c(d):= \sum_{k \in \Z^d} (1+|k|^2)^{2-n_*} <\infty.
\end{equation}

\noindent
\emph{Step 1} (Set-up). 
Let $u, w$ be as in the statement of the theorem.
Towards a contraposition
suppose that $\sup_{\overline{\cO}} (u-w)>0$. We fix  
a sufficiently small $\delta>0$ satisfying
$$
l :=\max_{(t,\mu)\in \overline{\cO}}\{(u-w)(t,\mu)- \delta(T-t)\}>0.
$$  
Set
$\bu(t,\mu):= u(t,\mu)-\delta(T-t)$.  Then,
$\bu$ is a continuous viscosity subsolution of
\begin{equation}
\label{eq:eta}
- \partial_t \bu(t,\mu)=
H(\mu,\pmu \bu(t,\mu))  -\delta.
\end{equation}

\noindent
\emph{Step 2} (Doubling the variables).
For $\eps>0$, set
$$
\Phi_{\eps}
(t,\mu,s,\nu):=\bu(t,\mu)-w(s,\nu)-\frac{1}{2\eps}\ 
\left(\rho_*^2(\mu,\nu) +(t-s)^2\right).
$$
As $\overline{\cO}$ is compact and $\bu, w$ are 
continuous, there exists $(t_\eps,s_\eps,\mu_\eps, \nu_\eps)
\in  \overline{\cO}\times \overline{\cO}$ satisfying  
$$
\Phi_\eps(t_\eps,\mu_\eps,s_\eps,\nu_\eps)
=\max_{ \overline{\cO}\times \overline{\cO}}\Phi_\eps\geq l >0.
$$
Set $M:= \max \bu$, $m:=\min v$, 
$\zeta_\eps:=\rho_*^2(\mu_\eps,\nu_\eps)
+(t_\eps-s_\eps)^2$,
so that
\begin{equation}
\label{eq:bound}
0\leq \zeta_\eps\leq 2 \eps\ (M+m-l).
\end{equation}

\noindent
\emph{Step 3} (Letting $\eps$ to zero).
Since $\overline{\cO}$ is compact, there is a subsequence
$\{(t_\eps,\mu_\eps,s_\eps,\nu_\eps)\}\subset 
\overline{\cO} \times \overline{\cO}$, 
denoted by $\eps$ again, and $(t^*,\mu^*,s^*,\nu^*) \in  \overline{\cO}\times \overline{\cO}$,
such that
$$
\mu_\eps \rightharpoonup\mu^*,\quad 
\nu_\eps\rightharpoonup \nu^*,\quad 
t_\eps \rightarrow  t^*,\quad 
s_\eps\rightarrow s^*,\qquad 
\text{as}\ \ \eps\downarrow 0.
$$
By \eqref{eq:bound} it is clear that $t^*=s^*$, 
and $\rho_*(\mu^*,\nu^*)=0$.  Then, by Lemma~\ref{lem:metric},
$\mu^*=\nu^*$.

If  $t^*$ were to be equal to $T$,
by the terminal condition 
\eqref{eq:initial}, we would have
$$
0<l \le \liminf_{\eps \downarrow 0}\,
\Phi_\eps(t_\eps,\mu_\eps,s_\eps,\nu_\eps)
\le \lim_{\eps \downarrow 0}\,
[\bar{u}(t_\eps,\mu_\eps) -w(s_\eps,\nu_\eps)]
= \bar{u}(T,\mu^*)- w(T,\mu^*) \le 0.
$$
Hence, $t^*<T$ and $t_\eps, s_\eps <T$ for all sufficiently
small $\eps>0$.
\vspace{2pt}

\noindent
\emph{Step 4} (Distance estimate).
Without loss of generality,
suppose that $w$ is Lipschitz, i.e,
$$
|w(t,\mu)-w(t,\nu)| \le \frac12 L_w\, \rho_*(\mu,\nu),
\qquad \mu,\nu \in \cP(\T^d), \, t \in [0,T].
$$  Then,
for each $\eps>0$,
\begin{align*}
 \bar{u}(t_\eps,\mu_\eps) - w(s_\eps,\nu_\eps)
 -\frac{1}{2 \eps}\, \zeta_\eps&=
 \Phi_\eps(t_\eps,\mu_\eps,s_\eps,\nu_\eps)
\ge  \Phi_\eps(t_\eps,\mu_\eps,s_\eps,\mu_\eps)\\
 &=
 \bar{u}(t_\eps,\mu_\eps) - w(s_\eps,\mu_\eps)
 -\frac{1}{2\eps}\, (t_\eps-s_\eps)^2.
 \end{align*}
 Therefore,
 $\rho^2_*(\mu_\eps,\nu_\eps)
 = \zeta_\eps - (t_\eps-s_\eps)^2 \le 2 \eps\,
 [w(s_\eps, \mu_\eps)-w(s_\eps, \nu_\eps)]
 \le 2\eps \, 
  L_w\, \rho_*(\mu_\eps,\nu_\eps)$.
Hence,
\begin{equation}
\label{eq:distest}
\rho_*(\mu_\eps,\nu_\eps) 
 \le \eps \, L_w\qquad
 \forall\, \eps>0.
 \end{equation}
 
\noindent
\emph{Step 5} (Viscosity property).
Set 
$$
\psi_\eps(t,\mu):=\frac{1}{2 \eps}
[\rho_*^2(\mu,\nu_\eps) +(t-s_\eps)^2],
\qquad 
\phi_\eps(s,\nu):=-\, \frac{1}{2 \eps}
[\rho_*^2(\mu_\eps,\nu) +(t_\eps-s)^2].
$$
By Lemma~\ref{lem:diff}, both 
$\partial_\mu\psi_\eps(t,\mu), \partial_\nu\phi_\eps(t,\mu) \in \cC^2(\T^d)$.  
Moreover, by the regularity
Assumption~\ref{a:regularity}, maps 
$(t,\mu) \mapsto H(\mu,\partial_\mu\psi_\eps(t,\mu))$, 
and $(t,\nu) \mapsto H(\nu,\partial_\mu\phi_\eps(t,\nu))$
are continuous.  Hence, $\psi_\eps$ and $\phi_\eps$
are smooth test functions. 
Set
$$
\kappa_\eps (x):=\partial_\mu \psi_\eps(t_\eps,\mu_\eps)(x)
=\partial_\mu \phi_\eps(s_\eps,\nu_\eps)(x)
= \frac{1}{\eps}\,
\sum_{k\in \Z^d} \,\frac{F_k(\mu_\eps-\nu_\eps)}
{(1+ |k|^2)^{n_*}}\ e^*_k(x),\qquad x \in {\T^d}.
$$

Also, $\bar{u}(t,\mu)- \psi_\eps(t,\mu)$
is maximized at $t_\eps, \mu_\eps$.  Since $t_\eps<T$, $\psi_\eps
\in \cC_s(\overline{\cO})$ and $\bar{u}$ is a viscosity 
subsolution of \eqref{eq:eta}, then
$$
- \frac{t_\eps-s_\eps}{\eps} \le H(\mu_\eps,\kappa_\eps) -\delta.
$$
By the viscosity property of $w$, a similar
argument implies that
$$
- \frac{t_\eps-s_\eps}{\eps} \ge H(\nu_\eps,\kappa_\eps).
$$
We subtract the above inequalities to arrive at
\begin{equation}
\label{eq:main}
0 < \delta \le  H(\mu_\eps,\kappa_\eps)
-  H(\nu_\eps,\kappa_\eps).
\end{equation}

\noindent
\emph{Step 6} (Estimation).
Since
$
H(\mu,\kappa_\eps) = \inf_{\alpha \in \cC_a} \left\{
\mu( \ell^\alpha(\cdot,\mu)+
\cM^{\alpha,\mu}[\kappa_\eps](\cdot))\right\},
$
$$
|H(\mu_\eps,\kappa_\eps) -  H(\nu_\eps,\kappa_\eps)|
\le 
\sup_{\alpha \in \cC_a}\cT^\alpha_\eps+
\sup_{\alpha \in \cC_a} \cI^\alpha_\eps+ 
\sup_{\alpha \in \cC_a}\cJ^\alpha_\eps,
$$
where
\begin{align*}
\cT^\alpha_\eps&:=  \left| \mu_\eps(\ell^\alpha(\cdot,\mu_\eps))
-  \nu_\eps(\ell^\alpha(\cdot,\nu_\eps)) \right|\\
\cI^\alpha_\eps&:=  \left| (\mu_\eps-\nu_\eps)(\cM^{\alpha,\mu_\eps}
[\kappa_\eps ](\cdot)) \right|\\
\cJ^\alpha_\eps&:=  \left| \nu_\eps(\cM^{\alpha,\mu_\eps}
[\kappa_\eps ](\cdot) -\cM^{\alpha,\nu_\eps}
[\kappa_\eps ](\cdot) ) \right|.
\end{align*}

\noindent
\emph{Step 7} (Estimating $\cT^\alpha_\eps$).
\vspace{2pt}
By the regularity Assumption~\ref{a:regularity}
and the estimate \eqref{eq:distest},
\begin{align*}
|\mu_\eps(\ell^\alpha(\cdot,\mu_\eps))
-  \nu_\eps(\ell^\alpha(\cdot,\nu_\eps))|&\leq
 |(\mu_\eps-\nu_\eps)(\ell^\alpha(\cdot,\mu_\eps))|
+ |\nu_\eps(\ell^\alpha(\cdot,\mu_\eps)
-  \ell^\alpha(\cdot,\nu_\eps))|\\
&\le \rho_*(\mu_\eps,\nu_\eps)\, \|\ell^\alpha(\cdot,\mu_\eps)\|_{\cC_*}
+\sup_{x \in {\T^d}} |\ell^\alpha(x,\mu_\eps)
-  \ell^\alpha(x,\nu_\eps)|\\
&\le 2c_a\, \rho_*(\mu_\eps,\nu_\eps)
\le 2c_a L_w\, \eps.
\end{align*}
Hence, we have $\lim_{\eps \downarrow 0}\ \sup_{\alpha \in \cC_a} \cT_{\eps}^\alpha=0$.
\vspace{4pt}

\noindent
\emph{Step 8} (Estimating $\cI_\eps$).
For $x \in {\T^d}$, $\mu \in \cP(\T^d)$, $\alpha \in \cC_a$, and $k \in \Z^d$
set 
$$
\beta^\alpha_{k}(x,\mu):=
\cM^{\alpha,\mu}[e^*_k](x)
 =-[k \cdot b^\alpha(x,\mu) + a_k^\alpha(x,\mu)] e^*_k(x),
$$
where
for $x\in \T^d, \mu \in \cP(\T^d), \alpha \in \cC_a, k \in \Z^d$,
\begin{equation}
\label{eq:ak}
a_k^\alpha(x,\mu):=
\frac{1}{2}
 \sum_{i,j=1}^d  \sum_{l=1}^{d'}\ \sigma_{il}(x,\mu,\alpha(x))
\sigma_{jl}(x,\mu,\alpha(x))  k_i k_j.
\end{equation}
Then, 
$$
\cM^{\alpha,\mu_\eps}[\kappa_\eps](x)
=\frac{1}{\eps}\sum_{k\in \Z^d}
\frac{1}{(1+|k|^2)^{n_*}}\ F_k(\mu_\eps-\nu_\eps)\
\beta^\alpha_{k}(x,\mu_\eps). 
$$
This in turn implies that
\begin{align*}
\cI^\alpha_\eps &\le 
\frac{1}{\eps}\sum_{k\in \Z^d}
\frac{1}{(1+|k|^2)^{n_*}}\ |F_k(\mu_\eps-\nu_\eps)|\
\left|(\mu_\eps-\nu_\eps)(\beta^\alpha_{k}(\cdot,\mu_\eps))\right|\\
&\le 
\frac{1}{\eps}\
(\sum_{k\in \Z^d}
\frac{|F_k(\mu_\eps-\nu_\eps)|^2}{(1+|k|^2)^{n_*}})^\frac12\
(\sum_{k\in \Z^d}
\frac{((\mu_\eps-\nu_\eps)(\beta^\alpha_{k}(\cdot,\mu_\eps)))^2}{(1+|k|^2)^{n_*}} 
)^\frac12\\
&\le 
\frac{\rho_*(\mu_\eps,\nu_\eps)}{\eps}\
(\sum_{k\in \Z^d} (1+|k|^2)^{2-n_*}\ \beta_{k,\eps}^2\ )^\frac12,
\end{align*}
where
$$
\beta_{k,\eps}:= (1+|k|^2)^{-1}\ \sup_{\alpha \in \cC_a} 
\left|(\mu_\eps-\nu_\eps)(\beta^\alpha_{k}(\cdot,\mu_\eps))
\right|,
\qquad  \in \Z.
$$

Again by Assumption~\ref{a:regularity},
$|\beta_{k,\eps}| \le c_a+c_a^2$, and $\beta^\alpha_{k,\eps}$ is
Lipschitz continuous with a Lipschitz constant
$c_k$ uniformly in $\alpha$.  Hence, by Kantorovich duality
$\beta_{k,\eps} \le c_k \hr_1(\mu_\eps,\nu_\eps)$.
As  $\mu_\eps-\nu_\eps$
converges weakly to zero, we conclude
that   $\beta_{k,\eps}$ also converges to zero
for every $k \in \Z$.
Also
$c(d)=\sum_{k=1}^{\infty}(1+|k|^2)^{2-n_*}$ is
finite by \eqref{eq:cd}, and we have argued that $|\beta_{k,\eps}|$ is
uniformly bounded. Hence, we may use
 dominated convergence to conclude that the sequence
$\sum_{k=1}^\infty (1+|k|^2)^{2-n_*}\beta_{k,\eps}^2$ 
converges to zero as $\eps \downarrow 0$. 
Then, by \eqref{eq:distest},
$$
\lim_{\eps \downarrow 0}\ \sup_{\alpha \in \cC_a} \cI^\alpha_\eps
\le \lim_{\eps \downarrow 0}\ 
L_w\
(\sum_{k=1}^\infty (1+|k|^2)^{2-n_*}\ \beta_{k,\eps}^2\ 
)^\frac12 =0.
$$

\noindent
\emph{Step 9} (Estimating $\cJ_\eps$). The definition of $\cJ^\alpha_\eps$
imply that
$$
\cJ^\alpha_\eps \le \sup_{x \in {\T^d}}\, \{
| \cM^{\alpha,\mu_\eps}[\kappa_\eps](x)
-\cM^{\alpha,\nu_\eps}[\kappa_\eps](x)|\}.
$$
Let $a^\alpha_k$ be as in \eqref{eq:ak},
and for $\alpha \in \cC_a$, $x \in {\T^d}$, $k \in \Z^d$, set
\begin{align*}
\gamma^\alpha_{k,\eps}(x)&:= \cM^{\alpha,\mu_\eps}[e^*_k](x)
-\cM^{\alpha,\nu_\eps}[e^*_k](x)\\
& =  k\cdot [b^\alpha(x,\nu_\eps) -b^\alpha(x,\mu_\eps)] e^*_k(x)
+ [a^\alpha_k(x,\nu_\eps) -a^\alpha_k(x,\mu_\eps)] e^*_k(x).
\end{align*}
By the regularity Assumption~\ref{a:regularity},
there exists $c_2$ such that
$$
\sup_{ x \in {\T^d}}|\gamma^\alpha_{k,\eps}(x)| 
\le c_2 (1+|k|^2) \hr_*(\mu_\eps,\nu_\eps),
\qquad \forall\alpha\in \cC_a, \  k \in \Z^d.
$$
 Hence, for every  $\balpha \in \cA$,
\begin{align*}
\cJ^\alpha_\eps
&\le \frac{1}{\eps}\ 
\sum_{k\in \Z^d}\ \frac{|F_k(\mu_\eps-
\nu_\eps)|}{(1+|k|^2)^{n_*}}\ \sup_{x \in {\T^d}}\,
|\gamma^{\alpha}_{k,\eps}(x)|\\
&\le 
\frac{c_2}{\eps}\,
\, (\sum_{k\in \Z^d}\ 
\frac{|F_k(\mu_\eps-\nu_\eps)|^2}{(1+|k|^2)^{n_*}})^\frac12\
(\sum_{k\in \Z^d}\ 
(1+|k|^2)^{2-n_*})^\frac12\,
\hr_*(\mu_\eps,\nu_\eps)\\
&\le  c_2 L_w\ c(d)\, \hr_*(\mu_\eps,\nu_\eps)
=:\hat{c}\ \hr_*(\mu_\eps,\nu_\eps),
\end{align*}
where $c(d)$ is as in \eqref{eq:cd}.
Therefore,
$\lim_{\eps \downarrow 0}\ \sup_{\alpha \in \cC_a}\cJ^\alpha_\eps 
\le \hat{c} \lim_{\eps \downarrow 0}\ \hr_*(\mu_\eps,\nu_\eps)=0$.
\vspace{4pt}

\noindent
\emph{Step 10} (Conclusion).   By \eqref{eq:main}
and above steps,
$ 0< \delta \le \lim_{\eps \downarrow 0}
\ [H(\mu_\eps,\kappa_\eps)
-  H(\nu_\eps,\kappa_\eps)] \le 0$.
This clear contradiction implies that 
$\max_{\overline{\cO}}\ (u-w) \le 0$.
\qed

\section{Lipschitz continuity}
\label{sec:Lip}
In this section, we prove Theorem~\ref{th:Lipschitz}.
\subsection{Regularity in space}
\label{ss:space}
We first prove
the continuous dependence of the 
solutions of the  McKean-Vlasov stochastic 
differential equation
\eqref{eq:mvsde} on its initial data.

\begin{proposition}
\label{pro:Lipsde}
Suppose that  the regularity Assumption~\ref{a:regularity}
holds.  Then,  there exists $\hc>0$ depending
on $T$ and the constant $c_a$ of Assumption~\ref{a:regularity}, such that
$$
\hr_*(\cL_u^{t,\mu,\balpha},\cL_u^{t,\nu,\balpha})
\le \hc\, \hr_*(\mu,\nu),\quad
\forall\ 0\le t\le u\le T, \ \mu, \nu \in \cP(\T^d),\ \balpha \in \cA.
$$
\end{proposition}
\begin{proof}  We complete the proof in several steps.
\vspace{4pt}

\noindent
\emph{Step 1} (Setting).
We fix $t \in [0,T]$,
$\mu, \nu \in \cP(\T^d)$, 
$\balpha \in \cA$, and set
$$
Y_u:= X_u^{t,\mu,\balpha}, \ \ \mu_u := \cL_u^{t,\mu,\balpha},
\qquad
Z_u:= X_u^{t,\nu,\balpha}, \ \ \nu_u := \cL_u^{t,\nu,\balpha},
\qquad u \in [t,T].
$$
By the definition of $\hr_*$, we need to prove the
following estimate for every $ u \in [t,T]$,
$$
(\mu_u-\nu_u)(\psi) \le \hc\, \hr_*(\mu,\nu)\,  \|\psi\|_{\cC_*},
\qquad \forall\, \psi \in \cC_*.
$$

\noindent
\emph{Step 2} (SDEs).
For $x \in \T^d$, let $Y^x, Z^x$ be the solutions of
the stochastic differential equations,
\begin{align*}
Y^x_u &= x + \int_t^u \left[b^{\alpha_s}(Y^x_s,\mu_s)\d s
+\sigma^{\alpha_s}(Y^x_s,\mu_s) \d B_s\right],\\
Z^x_u &= x + \int_t^u \left[b^{\alpha_s}(Z^x_s,\nu_s)\d s
+\sigma^{\alpha_s}(Z^x_s,\nu_s) \d B_s\right].
\end{align*}
Set $L_u^\mu(x):= \E[\psi(Y^x_u)]$, and 
$L_u^\nu(x):= \E[\psi(Z^x_u)]$.
Then, by conditioning, we have
$$
\mu_u(\psi) =\E[\psi(Y_u)] =\mu(L_u^\mu),
\qquad
\nu_u(\psi) =\E[\psi(Z_u)] =\nu(L_u^\nu).
$$
Therefore,
$$
(\mu_u-\nu_u)(\psi) = (\mu-\nu)(L^\mu_u) + \nu(L^\mu_u-L^\nu_u)
=:\cI_u(\psi) +\cJ_u(\psi).
$$

\noindent
\emph{Step 3} ($\cI_u$ estimate).
By the regularity Assumption~\ref{a:regularity},
there exists a constant $\hc_1$ satisfying
$$
\|b^{\alpha_u}(\cdot,\mu_u)\|_{\cC_*}
+\|\sigma^{\alpha_u}(\cdot,\mu_u)\|_{\cC_*}
\le \hc_1, \qquad \forall \ u \in [t,T].
$$
Hence, the map $x\in\T^d\to Y^x_u$ is 
$n_*$ times differentiable.  Therefore,
$L^\mu_u \in \cC_*$ and 
there exists a constant $\hc_2>0$ depending
only on $c_a$ of Assumption~\ref{a:regularity},
satisfying,
$$
\| L^\mu_u\|_{\cC_*} \le \hc_2\, \|\psi\|_{\cC_*}, \qquad \forall\ u \in [t,T], \mu \in \cP(\T^d).
$$
This implies that
$$
\cI_u(\psi) = (\mu-\nu)(L^\mu_u) \le  \hc_2\, \hr_*(\mu,\nu)\, \|\psi\|_{\cC_*} .
$$

\noindent
\emph{Step 4} ($\cJ_u$ estimate). By definitions, 
$\cJ_ \le \sup_{x} |L^\mu_u(x)-L^\nu_u(x)|$, and
$$
| L^\mu_u-L^\nu_u|
\le \E[|\psi(Y^x_u)-\psi(Z^x_u)|]
\le \E[|Y^x_u-Z^x_u|]\, \|\psi\|_1
\le  (\E[(Y_s^x-Z^x_s)^2])^\frac12\, \|\psi\|_*.
$$
For $x \in \T^d$, and  set 
$m^2_s(x):= \E[(Y_s^x-Z^x_s)^2]$.
We directly estimate that
$$
m^2_u(x) \le 2T \int_t^u \E[(b^{\alpha_s}(Y_s^x,\mu_s)-
b^{\alpha_s}(Z_s^x,\nu_s))^2]\d s
+2\int_t^u \E[|\sigma^{\alpha_s}(Y_s^x,\mu_s)-
\sigma^{\alpha_s}(Z_s^x,\nu_s)|^2]\d s.
$$
By the regularity Assumption~\ref{a:regularity}, 
$$
|b^{\alpha_s}(Y_s^x,\mu_s) -
b^{\alpha_s}(Z_s^x,\nu_s)|
\le c_a \left[ |Y_s^x-Z^x_s| +\hr_*(\mu_s,\nu_s)\right].
$$
Same estimate also holds for $|\sigma^{\alpha_s}(Y_s^x,\mu_s) -
\sigma^{\alpha_s}(Z_s^x,\nu_s)|$. Hence, there exists a constant
$\hc_3>0$, independent of $x$, satisfying, 
$m^2_u \le \hc_3 \int_t^u [m^2_s + \hr_*(\mu_s,\nu_s)^2] \d s$
for every $u \in [t,T]$.
By Gr\"onwall's inequality, there exists
$\hc_4>0$ satisfying 
$m^2_u \le  \hc_4^2 
\int_t^u \hr_*(\mu_s,\nu_s)^2 \d s$. Hence,
$$
\cJ_u \le  (\E[(Y_s^x-Z^x_s)^2])^\frac12\, \|\psi\|_* \le
\hc_4\, \left(
\int_t^u \hr_*(\mu_s,\nu_s)^2\ \d s\right)^\frac12
\, \|\psi\|_{\cC_*}, \qquad \forall u \in [t,T].
$$

\noindent
\emph{Step 5} (Conclusion). By the previous steps,
$$
(\mu_u-\nu_u)(\psi)  \le \left(  \hc_2\, \hr_*(\mu,\nu)\, +
 \hc_4\, \left(
\int_t^u \hr_*(\mu_s,\nu_s)^2\ \d s\right)^\frac12\right)\,
\|\psi\|_{\cC_*} ,
\qquad \forall \, \psi \in \cC_*.
$$
Since above holds for every $\psi \in \cC_*$,
the definition of $\hr_*$  implies that
$$
\hr_*(\mu_u,\nu_u) \le \hc_2\, \hr_*(\mu,\nu)\, +
 \hc_4\, \left(
\int_t^u \hr_*(\mu_s,\nu_s)^2\ \d s\right)^\frac12,
\qquad \forall \, u \in [t,T].
$$
Hence, 
$$
\hr_*(\mu_u,\nu_u)^2 \le 2 \hc_2^2\, \hr_*(\mu,\nu)^2\, +
2 \hc_4^2\int_t^u \hr_*(\mu_s,\nu_s)^2\d s,
\qquad \forall \, u \in [t,T].
$$
Again by Gr\"onwall, 
$\hr_*(\mu_u,\nu_u)^2
\le \hc^2\, \hr_*(\mu,\nu)^2$ for some $\hc>0$, for all $u \in [t,T]$.
\end{proof}

The following is an immediate consequence
of the above estimate.
\begin{lemma}
\label{lem:Lipschitz}
Under the regularity Assumption~\ref{a:regularity},
there exists
$L_1>0$ such that
$$
|J(t,\mu,\balpha)- J(t,\nu,\balpha)| \le 
L_1 \, \hr_*(\mu,\nu)\,
\qquad \forall\, \balpha \in \cA,\, 
\mu,\nu \in \cP(\T^d), \, t \in [0,T].
$$
Consequently, 
$$
|v(t,\mu)- v(t,\nu)| \le 
L_1 \, \hr_*(\mu,\nu)\,
\qquad \forall\, 
\mu,\nu \in \cP(\T^d), \, t \in [0,T].
$$
\end{lemma}
\begin{proof}
We fix $\balpha \in \cA$,
$\mu,\nu \in \cP(\T^d)$, $t \in [0,T]$,  and
use the same notation as in Proposition~\ref{pro:Lipsde}.
For $ u \in [t,T]$, the regularity 
Assumption~\ref{a:regularity} implies that
\begin{align*}
|\E[ \lau(Y_u,\mu_u) &-\lau(Z_u,\nu_u) ]|
\\&\le |\E[ \lau(Y_u,\mu_u) -\lau(Z_u,\mu_u) ]|
+ |\E[ \lau(Z_u,\mu_u) -\lau(Z_u,\nu_u) ]|\\
&\le |(\mu_u-\nu_u) (\lau(\cdot,\mu_u))|
+ c_a \,  \hr_*(\mu_u,\nu_u)\\
&\le \hr_*(\mu_u,\nu_u) \|\lau(\cdot,\mu_u)\|_{\cC_*}
 + c_a \,  \hr_*(\mu_u,\nu_u)\\
&\le 2 c_a\,  \hr_*(\mu_u,\nu_u) \le  2 c_a\, \hc \,  \hr_*(\mu,\nu).
\end{align*}
We now directly estimate using the above
to obtain the following inequalities,
\begin{align*}
|J(t,\mu,\balpha)- J(t,\nu,\balpha)|
&\leq \int_t^T |\E[ \lau(Y_u,\mu_u) -\lau(Z_u,\nu_u) ]|\ \d u\ +\
|\E[\varphi(\mu_T)
-\varphi(\nu_T)]|\\
&\le 2 c_a\, \hc \, (T-t)\,  \hr_*(\mu,\nu)
+ c_a \, \hr_*(\mu_T,\nu_T)\\
&\le c_a\, \hc \, ( 2(T-t) +1) \hr_*(\mu,\nu).
\end{align*}
As
$|v(t,\mu) -v(t,\nu)|
\le \sup_{\balpha \in \cA} |J(t,\mu,\balpha) -J(t,\nu,\balpha)|$,
the proof of the lemma is complete.

\end{proof}

\subsection{Time Regularity}
\label{ss:time}

\begin{proposition}
\label{pro:Liptime}
Suppose that  the regularity Assumption~\ref{a:regularity}
holds.  Then, there exists $L_2>0$ depending
on $T$ and the constant $c_a$ in Assumption~\ref{a:regularity}, such that
$$
|v(t,\mu) -v(\tau,\mu)|
\le L_2\, |t-\tau|^\frac12,\quad
\forall\ t, \tau\in [0,T], \ \mu \in \cP(\T^d).
$$
\end{proposition}
\begin{proof} 
Fix $0\le t \le \tau \le T$, $\mu \in \cP(\T^d)$, $\balpha \in \cA$,
and set $h:=\tau-t$.  With an arbitrary constant $a_*\in A$, we define
$$
\tilde{\alpha}_u(\cdot):=
\begin{cases}
\alpha_{u+h}(\cdot) & \text{if}\ u \in [t,T-h],\\
a_*& \text{if}\ u \in [T-h,T].
\end{cases}
$$
It is clear that $\tilde{\balpha} \in \cA$.
Set
$$
\tilde{\mu}_u:=\cL^{t,\mu,\tilde{\balpha}}_u,\quad u \in [t,T],
\quad\text{and}\quad
\mu_u:=\cL^{\tau,\mu,\balpha}_u,\quad u \in [\tau,T].
$$
Then, $\tilde{\mu}_u =\mu_{u+h}$
for every $u \in [t,T-h]$. In particular,
$$
\E[\ell^{\tilde{\alpha}_u}(X^{t,\mu,\tilde{\balpha}}_u)]=
\E[\lau(X^{\tau,\mu,\balpha}_{u+h})], \qquad \forall\, u \in [t,T-h].
$$

Since $\mu_T=\tilde{\mu}_{T-h}=\cL(X^{t,\mu,\tilde{\balpha}}_{T-h})$,
and $\tilde{\mu}_T=\cL(X^{t,\mu,\tilde{\balpha}}_{T})$,
$$
\hr_1(\tilde{\mu}_T,\mu_T)
\le \E[|X^{t,\mu,\tilde{\balpha}}_T-X^{t,\mu,\tilde{\balpha}}_{T-h}|]
\le \left(\E[(X^{t,\mu,\tilde{\balpha}}_T
-X^{t,\mu,\tilde{\balpha}}_{T-h})^2]\right)^\frac12.
$$
As $b, \sigma$ are bounded by $c_a$, 
there is $\tc_1>0$ satisfying,
$\hr_1(\tilde{\mu}_T,\mu_T)
\le \tc_1 \sqrt{h}$.  Therefore,
$$
|\varphi(\tilde{\mu}_T)-\varphi(\mu_T)|
\le c_a \hr_*(\tilde{\mu}_T,\mu_T)
\le c_a \hr_1(\tilde{\mu}_T,\mu_T)
\le \tc_1 \, c_a\, \sqrt{h}.
$$

Above estimate imply that for any 
$\balpha \in \cA$,
\begin{align*}
v(t,\mu) -J(\tau,\mu,\balpha)
&\le J(t,\mu,\tilde{\balpha}) -J(\tau,\mu,\balpha)\\
&= \int_{T-h}^T
\E[\ell^{\tilde{\alpha}_u}(X^{t,\mu,\tilde{\balpha}}_u)]\, \d u
+ \varphi(\tilde{\mu}_T)-\varphi(\mu_T)
\le c_a h + \tc_1 \, c_a\, \sqrt{h}.
\end{align*}
Hence,
$$
v(t,\mu) -v(\tau,\mu)
=\sup_{\balpha \in \cA} \left(v(t,\mu) -J(t,\mu,\balpha)\right)
\le c_a h + \tc_1 \, c_a\, \sqrt{h}.
$$

We prove the opposite inequality
by using the control
$$
\hat{\alpha}_u(\cdot):=
\begin{cases}
\alpha_{u-h}(\cdot) & \text{if}\ u \in [h,T],\\
a_*& \text{if}\ u \in [0,h].
\end{cases}
$$
Again $\hat{\balpha}\in \cA$, and we set
$$
\hat{\mu}_u:=\cL^{\tau,\mu,\hat{\balpha}}_u,\quad u \in [\tau,T],
\quad\text{and}\quad
\mu_u:=\cL^{t,\mu,\balpha}_u,\quad u \in [t,T].
$$
Then $\hat{\mu}_u=\mu_{u-h}$ for every $u\in[\tau,T]$ and $\hat{\mu}_T=\mu_{T-h}$. Following the above steps \emph{mutatis mutandis}, we obtain
the following inequality for any $\balpha\in\cA$,
\begin{align*}
v(\tau,\mu) -J(t,\mu,\alpha)
&\le J(\tau,\mu,\hat{\balpha}) -J(t,\mu,\alpha)\\
&= -\int_{t}^\tau
\E[\ell^{\tilde{\alpha}_u}(X^{t,\mu,\balpha}_u)]\, \d u
+ \varphi(\hat{\mu}_t)-\varphi(\mu_T)
\le c_a h + \tc_1 \, c_a\, \sqrt{h}.
\end{align*}
Hence,
$$
v(\tau,\mu) -v(t,\mu)
=\sup_{\balpha \in \cA} \left(v(\tau,\mu) -J(t,\mu,\balpha)\right)
\le c_a h + \tc_1 \, c_a\, \sqrt{h}.
$$

\end{proof}

\section{Dynamic Programming}
\label{sec:dpp}

In this section we prove Theorem~\ref{th:dpp}.  
For a general result but in a 
different setting, we refer the reader to \cite{DPT}.
\vspace{4pt}

\noindent\emph{Proof of Theorem~\ref{th:dpp}}.
We fix $(t,\mu)\in \overline{\cO}$, $\tau \in [t,T]$, and set
$$
Q(\balpha):= \int_t^\tau \E[\ell^{\alpha_s}(X^{t,\mu,\balpha}_s,
\cL^{t,\mu,\balpha}_s)]\d s 
+v(\tau,\cL^{t,\mu,\balpha}_\tau), \qquad \balpha \in \cA.
$$
Then, the dynamic programming principle can be stated as
$v(t,\mu)=\inf_{\balpha \in \cA} Q(\balpha)$.  
Recall that $v(t,\mu)=\inf_{\balpha \in \cA} J(t,\mu,\balpha)$. 
For any 
$\balpha \in\cA$, and $s \in [\tau,T]$, Markov property implies that
$X^{t,\mu,\balpha}_s = X^{\tau,\cL^{t,\mu,\balpha}_\tau,\balpha}_s$,
and consequently
$\cL^{t,\mu,\balpha}_s = \cL^{\tau,\cL^{t,\mu,\balpha}_\tau,\balpha}_s$.
Hence,
\begin{align*}
\int_\tau^T \E[\ell^{\alpha_s}(X^{t,\mu,\balpha}_s,\cL^{t,\mu,\balpha}_s)]\ \d s 
&+\varphi(\cL^{t,\mu,\balpha}_T) \\
&= \int_\tau^T \E[\ell^{\alpha_s}(X^{\tau,\cL^{t,\mu,\balpha}_\tau,\balpha}_s,
\cL^{\tau,\cL^{t,\mu,\balpha}_\tau,\balpha}_s)]\ \d s 
+\varphi(\cL^{\tau,\cL^{t,\mu,\balpha}_\tau,\balpha}_T) \\
&=J(\tau, \cL^{\tau,\cL^{t,\mu,\balpha}_\tau,\balpha}_\tau,\balpha)
\ge v(\tau, \cL^{\tau,\cL^{t,\mu,\balpha}_\tau,\balpha}_\tau).
\end{align*}
This implies that 
\begin{align*}
J(t,\mu,\balpha) &=
\int_t^\tau \E[\ell^{\alpha_s}(X^{t,\mu,\balpha}_s,\cL^{t,\mu,\balpha}_s)]\ \d s+
\left(\int_\tau^T \E[\ell^{\alpha_s}(X^{t,\mu,\balpha}_s,\cL^{t,\mu,\balpha}_s))]\ \d s 
+\varphi(\cL^{t,\mu,\balpha}_T) \right) \\
&\ge \int_t^\tau \E[\ell^{\alpha_s}(X^{t,\mu,\balpha}_s,\cL^{t,\mu,\balpha}_s))]\ \d s
 +v(\tau, \cL^{\tau,\cL^{t,\mu,\balpha}_\tau,\balpha}_\tau)= Q(\balpha).
 \end{align*}
 Therefore,
 $ v(t,\mu)=\inf_{\balpha \in \cA} J(t,\mu,\balpha) 
 \ge \inf_{\balpha \in \cA} Q(\balpha)$.

To prove the opposite inequality, we 
fix $\eps>0$, and  set $\delta := \eps/(4 L_1)$.
By Lemma~\ref{lem:Lipschitz}, 
whenever $\hr_*(\nu,\eta)\le \delta$,  we have
$|J(\tau,\nu,\balpha)- J(\tau,\eta,\balpha)| 
\le \eps/4$, for every $\balpha \in \cA$,
and also $|v(\tau,\nu)- v(\tau,\eta)| \le \eps/4$.
Consider a covering
of $\cP(\T^d)$ given by
$$
\cB(\nu):= \{ \eta \in  \cP(\T^d)\ :\
\hr_*(\nu,\eta) <\delta\ \},
\qquad \nu \in \cP(\T^d).
$$
It is clear that each $\cB(\nu)$ is an open set
as  $\hr_*$ is continuous with respect to
the weak$^*$ topology. Then, since  $\cP(\T^d)$ is
weak$^*$  compact, 
there exits 
$\{\nu_j\}_{j=1,\ldots,n}\subset \cP(\T^d)$
such that
$\cP(\T^d) = \cup_{j=1}^{n}\ \cB(\nu_j)$.
Set $\cB_1:= \cB(\nu_1)$, and
and recursively define
$$
\cB_{j+1} := \cB(\nu_{j+1}) \setminus \cup_{l=1}^j \cB_l,
\qquad j=1,\ldots, n-1,
$$ 
so that $\{\cB_j\}_{j=1,\ldots,n}$ forms a disjoint
covering of $\cP(\T^d)$. Moreover, 
for any $\nu \in \cB_j \subset \cB(\nu_j)$,
$\hr_*(\nu,\nu_j)\le \delta$, and therefore,
$$
|v(\tau,\nu)-v(\tau,\nu_j)| \le \frac{\eps}{4},
\quad\text{and}\quad
|J(\tau,\nu,\balpha)-J(\tau,\nu_j,\balpha)| \le \frac{\eps}{4},
\quad\forall\ \balpha \in \cA.
$$

For each $j$, choose $\balpha^j \in \cA$ so that
$J(\tau,\nu_j,\balpha^j) \le v(\tau,\nu_j)+\frac{\eps}{4}$.
Then,
\begin{equation}
\label{eq:vq}
J(\tau,\nu,\balpha^j) \le 
J(\tau,\nu_j,\balpha^j) +\frac{\eps}{4}
\le v(\tau,\nu_j)+\frac{\eps}{2}
\le v(\tau,\nu)+\frac{3\eps}{4},
\qquad \forall\ \nu \in \cB_j.
\end{equation}
We  choose $\balpha^*\in \cA$ satisfying
$Q(\balpha^*) \le \inf_{\balpha \in \cA} Q(\balpha) + \frac{\eps}{4}$,
and define a control process $\balpha^\eps$ by,
$$
\alpha^\eps_u(x)=
\begin{cases} 
\alpha^*_u(x), &\quad \text{if}\  u \in [t,\tau), \\
& \\
\sum_{j=1}^n \alpha^j_u(x) \chi_{\cB_j}
(\cL^{t,\mu,\balpha^*}_\tau),
 &\quad \text{if}\ u \in[\tau,T],	
\end{cases}
\qquad x \in \T^d.
$$

As $\balpha^*$ and $\balpha^\eps$ agree on $[t,\tau]$,
we have $\cL^{t,\mu,\balpha^*}_u=\cL^{t,\mu,\balpha^\eps}_u$ 
for all $u \in[t,\tau]$.  Hence,
$$
 \inf_{\balpha \in \cA} Q(\balpha)+\frac{\eps}{4} \ge 
Q(\balpha^*) 
=  \int_t^\tau \E[\ell^{\alpha^*_s}(X^{t,\mu,\balpha^*}_s,\cL^{t,\mu,\balpha^*}_s)] \d s 
+v(\tau,\cL^{t,\mu,\balpha^*}_\tau) = Q(\balpha^\eps).
$$
Moreover, by the definition of $\balpha^\eps$
and \eqref{eq:vq},
\begin{align*}
 v(\tau,\cL^{t,\mu,\balpha^*}_\tau)&=
 \sum_{j=1}^n v(\tau,\cL^{t,\mu,\balpha^*}_\tau)
\chi_{\cB_j}(\cL^{t,\mu,\balpha^*}_\tau ) 
\ge
\sum_{j=1}^n J(\tau,\cL^{t,\mu,\balpha^*}_\tau, \balpha^j)
\chi_{\cB_j}(\cL^{t,\mu,\balpha^*}_\tau ) -\frac{3\eps}{4}\\
&= J(\tau,\cL^{t,\mu,\balpha^*}_\tau, \balpha^\eps)-\frac{3\eps}{4}.
\end{align*}
Hence,
\begin{align*}
 \inf_{\balpha \in \cA} Q(\balpha) +\eps &\ge
 Q(\balpha^\eps) +\frac{3 \eps}{4}
 = \int_t^\tau \E[\ell^{\alpha^*_s}(X^{t,\mu,\balpha^*}_s,\cL^{t,\mu,\balpha^*}_s)] \d s 
+\left(v(\tau,\cL^{t,\mu,\balpha^*}_\tau)+\frac{3\eps}{4}\right)\\
&\ge \int_t^\tau  \E[\ell^{\alpha^*_s}(X^{t,\mu,\balpha^*}_s,\cL^{t,\mu,\balpha^*}_s)] \d s  
+J(\tau,\cL^{t,\mu,\balpha^\eps}_\tau, \balpha^\eps)\\
&=J(t,\mu,\balpha^\eps)
\ge v(t,\mu).
\end{align*}
\qed

\section{Viscosity property}
\label{sec:viscosity}

In this section, we prove the viscosity property 
of the value function.  Although the 
below  proof follows the standard 
one very closely, 
we provide it for completeness.

The following version of the It\^o's formula along flows of measures
follows from Proposition 5.102 of \cite{CD}.  Recall that $X^{t,\mu,\balpha}$
is the solution of \eqref{eq:mvsde}, 
$\cL^{t,\mu,\balpha}_u= \cL(X^{t,\mu,\balpha}_u)$,
and the operator  $\cM^{a,\mu}$ is defined in subsection~\ref{ss:dpp}.

\begin{lemma}
\label{lem:ito}
For every $\psi\in \cC_s(\T^d)$, $(t,\mu)\in \overline{\cO}$,
$u \in [t,T]$, and $\balpha\in \cA$, 
$$
\psi(u,\cL^{t,\mu,\balpha}_u)=\psi(t,\mu)+\int_t^u
\left( \partial_t\psi(s,\cL^{t,\mu,\balpha}_s)
+\E[\cM^{\alpha_s,\cL^{t,\mu,\balpha}_s}[\partial_\mu \psi(s,\cL^{t,\mu,\balpha}_s)](X^{t,\mu,\balpha}_s)]
\right)\d s.
$$
\end{lemma}

\subsection{Subsolution}
\label{ss:sub}

Suppose that for $(t_0,\mu_0)\in [0,T)\times \cP(\T^d)$ 
and  test function $\psi\in \cC_s(\overline{\cO})$,
$$
0=(v-\psi)(t_0,\mu_0)=\max_{\overline{\cO}} (v-\psi).
$$
For $\alpha \in \cC_a$, set
$$
k^\alpha(t,x,\mu):= 
\ell(x,\mu,\alpha(x))+\cM^{\alpha,\mu}[\partial_\mu \psi(t,\mu)])](x), 
\qquad t \in [0,T],\ x \in \T^d,\, \mu \in \cP(\T^d).
$$
As
$H(\mu_0,\partial_\mu \psi(t_0,\mu_0))= \inf_{\alpha \in \cC_a}
\mu_0(k^\alpha(t_0,\cdot,\mu_0))$,
for any $\eps>0$ there is $\alpha^* \in \cC_a$ satisfying,
$$
\mu_0(k^{\alpha^*}(t_0,\cdot,\mu_0)) 
\le H(\mu_0,\partial_\mu \psi(t_0,\mu_0))+\eps.
$$

Set $\alpha^*_u\equiv \alpha_*$ and 
let $X^*_u:= X^{t_0,\mu_0,\balpha^*}_u$ and 
$\mu^*_u:= \cL^{t_0,\mu_0,\balpha^*}_u$
for $ u \in [t_0,T]$.
Since $v \le \psi$,
dynamic programming principle Theorem~\ref{th:dpp}
with $\tau=t_0+h\leq T$ implies that
$$
v(t_0,\mu_0) \leq  \int_{t_0}^{t_0+h}
\E[\ell(X^*_s,\mu^*_s,\alpha^*(X^*_s)] \d s
+\psi(t_0+h,\mu^*_{t_0+h}).
$$
By Lemma~\ref{lem:ito},
$$
\psi(t_0+h,\mu^*_{t_0+h})
= \psi(t_0,\mu_0)+\int_{t_0}^{t_0+h}
\left( \partial_t\psi(s,\mu^*_s)
+\E[\cM^{\alpha^*,\mu^*_s}
[\partial_\mu \psi(s,\mu^*_s)] (X^*_s)] \right) \d s.
$$

Since $\psi(t_0,\mu_0)=v(t_0,\mu_0)$, above inequalities
imply that
\begin{equation}
\label{eq:sub1}
0 \le  \frac{1}{h}\int_{t_0}^{t_0+h}
\left( \partial_t\psi(s,\mu^*_s)
+\E[ k^{\alpha^*}(s, X^*_s,\mu^*_s)]\right)\d s.
\end{equation}
We now let $h$ tend to zero to arrive at
the following inequality,
$$
-  \partial_t\psi(t_0,\mu_{0})  \le
\E[ k^{\alpha^*}(t_0, X_{t_0},\mu_{t_0})]
= \mu_0(k^{\alpha^*}(t_0,\cdot,\mu_0))
 \le H(\mu_0,\partial_\mu \psi(t_0,\mu_0))+\eps.
$$

\subsection{Supersolution}
\label{ss:super}

Suppose that for $(t_0,\mu_0)\in[0,T)\times\cP(\T^d)$ 
and a test function $\psi\in\cC_s(\overline{\cO})$, 
$$
0=(v-\psi)(t_0,\mu_0)=\min_{\overline{\cO}}(v-\psi).
$$
We may assume that the minimum is strict.
Towards a counterposition, suppose that
$$
- \partial_t\psi(t_0,\mu_0)< H(\mu_0,\partial_\mu\psi(t_0,\mu_0))
= \inf_{\alpha \in \cC_a} \left\{
\mu_0(k^\alpha,\cdot,\mu_0))\right\},
$$
where $k^\alpha(t,x,\mu)= \ell^\alpha(x,\mu)+
\cM^{\alpha,\mu}[\partial_\mu\psi(t,\mu)](x)$
is as in the previous subsection.
By  Definition~\ref{def:smooth} of test functions $\cC_s(\cO)$,
the map $(t,\mu) \in \overline{\cO}
\mapsto H(\mu, \partial_\mu \psi(t,\mu))$ is continuous.  
Therefore, there exists $\delta>0$
and a neighborhood $\cB\subseteq \overline{\cO}$ of $(t_0,\mu_0)$ 
such that 
$$
- \partial_t\psi(s,\mu) +\delta
\le  H(\mu,\partial_\mu\psi(t,\mu))
= \inf_{\alpha \in \cC_a} \left\{
\mu(k^\alpha(t,\cdot,\mu))\right\},\qquad \forall\, (t,\mu)\in\cB.
$$

For $\balpha \in \cA$, 
set $X^\balpha_s:=X^{t_0,\mu_0,\balpha}_s$,
$\mu^\balpha_s:=\cL^{t_0,\mu_0,\balpha}_s$,
and consider the (deterministic) time
$$
\tau^\balpha:=\inf\{s \in [t_0,T]: (s,\mu^{\balpha}_s) \notin \cB\},
$$
so that for every $s \in [t_0,\tau^\balpha)$, $(s,\mu^\balpha_s )\in \cB$, 
and consequently
$$
 \mu^{\balpha}_s( k^{\alpha_s}(s,\cdot,\mu^\balpha_s))
\ge H(\mu^{\balpha}_s,\partial_\mu\psi(s,\mu^{\balpha}_s))
\ge - \partial_t\psi(s,\mu^{\balpha}_s)+\delta.
$$
As $ \E[ k^{\alpha_s}(s,X_s^\balpha,\mu^\balpha_s)]
=\mu^{\balpha}_s( k^{\alpha_s}(s,\cdot,\mu^\balpha_s))$,
$$
\int_{t_0}^{\tau^\balpha} \left(\E[ k^{\alpha_s}(X_s^\balpha,\mu^\balpha_s)]
+ \partial_t\psi(s,\mu^{\balpha}_s)\, \right)\, \d s
\ge \delta (\tau^\balpha-t_0).
$$
Then, by Lemma~\ref{lem:ito}, we obtain the following inequality,
\begin{align*}
\psi(\tau^\balpha,\mu^\balpha_{\tau^\balpha})
&= \psi(t_0,\mu_0)
+\int_{t_0}^{\tau^\balpha} ( \partial_t\psi(s,\mu^{\balpha}_s)
+ \E[ \cM^{\alpha_s,\mu^\balpha_s}[\partial_\mu \psi(s,\mu^\balpha_s)](X_s^\balpha)\ ]) \d s\\
&= \psi(t_0,\mu_0)
+\int_{t_0}^{\tau^\balpha} \left( \partial_t\psi(s,\mu^{\balpha}_s) 
+\E[ k^{\alpha_s}(s,X_s^\balpha,\mu^\balpha_s)]
- \E[ \ell^{\alpha_s}(X_s^\balpha,\mu^\balpha_s)]\right) \d s\\
&\ge \psi(t_0,\mu_0)
-\int_{t_0}^{\tau^\balpha} \E[ \ell^{\alpha_s}(X_s^\balpha,\mu^\balpha_s)]\ \d s
+\delta (\tau^\balpha-t_0).
\end{align*}

Since $v \ge \psi$ and
$\psi(t_0,\mu_0) =v(t_0,\mu_0)$, above implies that
$$
v(t_0,\mu_0)
\le \int_{t_0}^{\tau^\balpha} 
\E[ \ell^{\alpha_s}(X_s^\balpha,\mu^\balpha_s)]\ \d s
+ v(\tau^\balpha,\mu^\balpha_{\tau^\balpha})-g(\balpha),
\qquad \forall \ \balpha \in \cA.
$$
where $g(\balpha):=\delta(\tau^\balpha-t_0)
+(v(\tau^\balpha,\mu^\balpha_{\tau^\balpha})
-\psi(\tau^\balpha,\mu^\balpha_{\tau^\balpha}))$.
We now claim that
$$
\delta_0:= \inf_{\balpha \in \cA}  g(\balpha) >0.
$$
Indeed, since $v \ge \psi$, if  $\tau^\balpha=T$, 
then $g(\balpha)\ge \delta(T-t_0)$.
On the other hand if $\tau^\balpha<T$, then 
$(\tau^\balpha,\mu^\balpha_{\tau^\balpha})\in \partial \cB$.
As $\cB$ is compact and $(t_0,\mu_0) \notin \partial \cB$
is the strict minimizer of $v-\psi$, we have 
$$
(v-\psi)(\tau^\balpha,\mu^\balpha_{\tau^\balpha})
\ge \inf_{(t,\mu) \in \partial \cB} (v-\psi)(t,\mu) >0.
$$
Hence, $\delta_0>0$ and the
above inequalities imply that
for \emph{every} $\balpha \in \cA$,
$$
v(t_0,\mu_0)
\le \int_{t_0}^{\tau^\balpha} \E[ \ell^{\alpha_s}(X_s^\balpha,\mu^\balpha_s)]\\ \d s
+ v(\tau^\balpha,\mu^\balpha_{\tau^\balpha})-\delta_0.
$$
This contradiction
to dynamic programming  implies that 
$-\psi_t(t_0,\mu_0)\ge H(\mu_0,
\partial_\mu \psi(t_0,\mu_0))$.

\qed
\vskip 0.2in

\bibliographystyle{abbrvnat}
\bibliography{wasserstein}

\end{document}